# Double Integrals for Euler's Constant and $\ln\frac{4}{\pi}$ and an Analog of Hadjicostas's Formula

## Jonathan Sondow

**1. INTRODUCTION.** Euler's constant, $\gamma$, is defined as the limit

$$\gamma = \lim_{N\to\infty}\left(1 + \frac{1}{2} + \frac{1}{3} + \cdots + \frac{1}{N} - \ln N\right). \tag{1}$$

(For a proof that the limit exists, see [**3**, sec. 10.1], [**7**, Chap. 2] or [**11**, sec. 9.1].) In this note we prove the formulas

$$\gamma = \sum_{n=1}^{\infty}\left(\frac{1}{n} - \ln\frac{n+1}{n}\right) = \iint_{[0,1]^2} \frac{1-x}{(1-xy)(-\ln xy)}\,dx\,dy, \tag{2}$$

$$\ln\frac{4}{\pi} = \sum_{n=1}^{\infty}(-1)^{n-1}\left(\frac{1}{n} - \ln\frac{n+1}{n}\right) = \iint_{[0,1]^2} \frac{1-x}{(1+xy)(-\ln xy)}\,dx\,dy. \tag{3}$$

In view of series (2), which is due to Euler, series (3) reveals $\ln(4/\pi)$ to be an "alternating Euler constant."

The constants $\gamma$ and $\ln(4/\pi)$ are related by Euler's formula

$$\gamma = \ln\frac{4}{\pi} + 2\sum_{n=2}^{\infty}(-1)^n \frac{\zeta(n)}{2^n n}, \tag{4}$$

where

$$\zeta(s) = 1 + \frac{1}{2^s} + \frac{1}{3^s} + \cdots \quad (\Re(s) > 1) \tag{5}$$

is the Riemann zeta function. Relation (4), a special case of a formula for the gamma function [**1**, eq. 6.1.33], can be proved using series (2), (3) and (5).

We constructed the double integral (2) for $\gamma$ by analogy with Beukers' integrals [**2**]

$$\zeta(2) = \iint_{[0,1]^2}\frac{1}{1-xy}\,dx\,dy, \quad \zeta(3) = \frac{1}{2}\iint_{[0,1]^2}\frac{-\ln xy}{1-xy}\,dx\,dy. \tag{6}, (7)$$

Since going "up" from $\zeta(2)$ to $\zeta(3)$ involves multiplying the integrand of (6) by the derivative

$$-\ln xy = -\frac{\partial}{\partial t}(xy)^t\bigg|_{t=0},$$

in order to go "down" from $\zeta(2)$ to $\gamma$ (which by comparing (1) and (5) one may think of as "$\zeta(1)$") we multiplied by the integral

$$-\frac{1}{\ln xy} = \int_0^\infty (xy)^t\, dt. \tag{8}$$

Because of the asymmetry in definition (1) (and the preference for a convergent improper integral), we also multiplied by $(1-x)$. It then turned out that the resulting integral (2) does indeed represent $\gamma$.

We noticed (3) in a different way. We first made series (2) alternate and correspondingly changed a sign in the denominator of integral (2). We then found the value, $\ln(4/\pi)$, of the modified series and integral. Related formulas for $\ln(\pi/2)$ are given in [**17**].

Beukers used "damped" versions of integrals (6) and (7) to simplify Apéry's famous proof of the irrationality of $\zeta(2)$ and $\zeta(3)$ (see [**2**], [**4**, sec. 11.3], [**12**] and [**19**]). "Damping" integrals (2) and (3), however, only leads to *criteria* for irrationality of $\gamma$ [**16**] and $\ln(4/\pi)$; their arithmetic nature remains undetermined.

We prove the series formulas (2) and (3) in Section 2, and the integral ones in Section 3. We state a recent generalization of the integral formula (2) in Section 4, and prove an analogous generalization of (3) in the final section. A special case is the evaluation

$$\iint_{[0,1]^2} \frac{1-x}{(1+xy)(\ln xy)^2}\, dxdy = \ln\frac{\pi^{1/2}A^6}{2^{7/6}e}, \tag{9}$$

where *A* is the Glaisher-Kinkelin constant, defined as the limit (see [**8**, sec. 2.15])

$$A = \lim_{n\to\infty} \frac{1^1 2^2 3^3 \cdots n^n}{n^{\frac{1}{2}n^2+\frac{1}{2}n+\frac{1}{12}} e^{-\frac{1}{4}n^2}}. \tag{10}$$

The constant *A* plays the same role in the approximation (10) to $1^1 2^2 3^3 \cdots n^n$ as the constant $\sqrt{2\pi}$ plays in Stirling's approximation to $n!$, which can be written

$$\sqrt{2\pi} = \lim_{n\to\infty} \frac{1\cdot 2\cdot 3\cdots n}{n^{n+1/2}e^{-n}}.$$

Euler's constant plays a similar role in the approximation (1) to $1+\frac{1}{2}+\frac{1}{3}+\cdots+\frac{1}{n}$ by $\ln n$.

**2. SUMMING THE SERIES.** To see that Euler's series (2) sums to $\gamma$, write the series as the limit

$$\lim_{N\to\infty} \sum_{n=1}^{N} \left(\frac{1}{n} - \ln\frac{n+1}{n}\right) = \lim_{N\to\infty} \left(\sum_{n=1}^{N} \frac{1}{n} - \ln N - \ln\frac{N+1}{N}\right)$$

and use definition (1).

To sum series (3), write it as the difference between the alternating harmonic series,

$$\sum_{n=1}^{\infty} (-1)^{n-1} \frac{1}{n} = 1 - \frac{1}{2} + \frac{1}{3} - \cdots = \ln 2, \tag{11}$$

and the logarithm of Wallis's product,

$$\prod_{n=1}^{\infty} \left(\frac{n+1}{n}\right)^{(-1)^{n-1}} = \frac{2}{1}\frac{2}{3}\frac{4}{3}\frac{4}{5}\frac{6}{5}\frac{6}{7}\cdots = \frac{\pi}{2},$$

obtaining $\ln 2 - \ln(\pi/2) = \ln(4/\pi)$.

**3. EVALUATING THE INTEGRALS.** Denote integral (2) by $I$ and expand $1/(1-xy)$ in a geometric series:

$$I = \iint_{[0,1]^2} \sum_{n=0}^{\infty} \frac{1-x}{-\ln xy} (xy)^n \, dxdy = \sum_{n=0}^{\infty} \iint_{[0,1]^2} (1-x) \frac{(xy)^n}{-\ln xy} \, dxdy.$$

(To justify interchanging summation and integration, replace the series on the left by a finite sum plus remainder and estimate the remainder—see [**16**, sec. 2].) Substitute the generalization of (8)

$$-\frac{(xy)^n}{\ln xy} = \int_n^{\infty} (xy)^t \, dt$$

and reverse the order of integration (permitted because the integrand is non-negative), obtaining

$$I = \sum_{n=0}^{\infty} \iint_{[0,1]^2} \int_n^{\infty} (1-x)(xy)^t \, dt\, dx\, dy = \sum_{n=0}^{\infty} \int_n^{\infty} \iint_{[0,1]^2} (x^t - x^{t+1}) y^t \, dx\, dy\, dt.$$

Replace $n$ by $n-1$ and integrate, first with respect to $x$ and $y$, then with respect to $t$, to arrive at series (2) for $\gamma$. (For a different proof, see the remark at the end of [**18**].)

The evaluation of integral (3) is similar. Since the expansion of $1/(1+xy)$ is an alternating series, the result is series (3) for $\ln(4/\pi)$.

**4. HADJICOSTAS'S FORMULA.** After seeing the integral formula (2) for $\gamma$, Hadjicostas conjectured a generalization [**9**], which Chapman then proved [**5**]:

$$\iint_{[0,1]^2} \frac{1-x}{1-xy}(-\ln xy)^s \, dx \, dy = \Gamma(s+2)\left[\zeta(s+2) - \frac{1}{s+1}\right] \quad (\Re(s) > -2). \tag{12}$$

Here $\Gamma$ is defined by Euler's integral [**6**],

$$\Gamma(s) = \int_0^\infty e^{-t} t^{s-1} dt \quad (\Re(s) > 0), \tag{13}$$

and $\zeta$ is the analytic continuation of (5) (see, e.g., [**13**] or [**17**]).

Taking the limit as $s \to -1$ in (12) and using the fact that (see, e.g., [**14**])

$$\lim_{s \to 1}\left(\zeta(s) - \frac{1}{s-1}\right) = \gamma, \tag{14}$$

we obtain (2). Conversely, (2) and (12) imply (14).

For $s = 0$ and 1, formula (12) is related to the integral representations (6) and (7) of $\zeta(2)$ and $\zeta(3)$, respectively.

**5. AN ANALOG.** Hadjicostas asked if there is an analog of (12) for the integral formula (3) for $\ln(4/\pi)$ [**10**]. The following result fills the bill and, as a bonus, yields (9).

**Theorem 1.** *For complex s with $\Re(s) > -3$, we have*

$$\iint_{[0,1]^2} \frac{1-x}{1+xy}(-\ln xy)^s \, dx \, dy = \Gamma(s+2)\left[\zeta^*(s+2) + \frac{1-2\zeta^*(s+1)}{s+1}\right], \tag{15}$$

*where $\zeta^*(s)$, the alternating zeta function, is the analytic continuation of the convergent Dirichlet series*

$$\zeta^*(s) = 1 - \frac{1}{2^s} + \frac{1}{3^s} - \cdots \quad (\Re(s) > 0).$$

Here $\Gamma(s)$ is analytically continued for $\Re(s) > -1$ and $s \neq 0$ by the functional equation (established by integrating by parts in (13))

$$\Gamma(s) = \frac{\Gamma(s+1)}{s}, \tag{16}$$

and $\zeta^*(s)$ is an entire function related to $\zeta(s)$ by the product (see, e.g., [**13**] or [**15**])

$$\zeta^*(s) = (1 - 2^{1-s})\zeta(s). \tag{17}$$

*Proof of Theorem 1.* The integral in (15) defines a holomorphic function $I(s)$ on the half-plane $\Re(s) > -3$. We prove (15) for $\Re(s) > 0$ and the theorem then follows by analytic continuation.

We make the change of variables $u = xy$, $v = 1 - x$ and integrate with respect to $v$. We then substitute $t = -\ln u$, obtaining

$$I(s) = \int_0^\infty \left( \frac{t^{s+1} - 2t^s}{e^t + 1} + e^{-t} t^s \right) dt.$$

Using the evaluation

$$\int_0^\infty \frac{t^{s-1}}{e^t + 1} dt = \Gamma(s)\zeta^*(s) \tag{18}$$

(to prove (18), write the denominator as $e^t(1 + e^{-t})$ and expand $1/(1 + e^{-t})$ in a geometric series, then integrate termwise via (13)), we get

$$I(s) = \Gamma(s+2)\zeta^*(s+2) + \Gamma(s+1)\left[1 - 2\zeta^*(s+1)\right].$$

Now (16) (with $s+1$ in place of $s$) gives (15) and the theorem follows.   ●

In order to show that Theorem 1 implies the integral formulas (3) and (9), we need the following values of the alternating zeta function and its derivative: $\zeta^*(1) = \ln 2$ (from (11)), $\zeta^*(0) = 1/2$, $\zeta^*(-1) = 1/4$, and $\zeta^{*\prime}(0) = \frac{1}{2}\ln\frac{\pi}{2}$ (see [**13**]). It follows that $I(-1) = \zeta^*(1) - 2\zeta^{*\prime}(0)$, which gives (3). Using (16) (with $s+2$ in place of $s$), we deduce that $I(-2) = \zeta^{*\prime}(0) + 2\zeta^*(-1) - 1 + 2\zeta^{*\prime}(-1)$. Employing (17) and the formula $\ln A = \frac{1}{12} - \zeta'(-1)$ for the Glaisher-Kinkelin constant (see [**8**, sec. 2.15]), we arrive at (9).

*209 West 97th Street, New York, NY 10025*
*jsondow@alumni.princeton.edu*